\documentclass[11pt]{article}

\usepackage[latin1]{inputenc}
\usepackage{url, amsmath,enumerate,fancyhdr,amssymb, amsthm, 
pstricks, epsf, lscape, ifpdf, graphicx, float}
\usepackage{lineno}  

\evensidemargin \oddsidemargin
\advance \topmargin by -\headheight
\advance \topmargin by -\headsep
\newcommand{\compilepdf}[1]{#1} 

\mathchardef\mhyphen="2D

\newtheorem{Definition}{Definition}
\newtheorem{Example}{Example}
\newtheorem{Proposition}{Proposition}
\newtheorem{Lemma}{Lemma}
\newtheorem{Theorem}{Theorem}
\newtheorem{Corollary}{Corollary}
\newtheorem{Remark}{Remark}
\newtheorem{Assumption}{Assumption}
\newtheorem{Conjecture}{Conjecture}

\newcommand{\lin}{\operatorname{lin}}
\newcommand{\cl}{\operatorname{cl}}
\newcommand{\ri}{\operatorname{ri}}
\newcommand{\dir}{\operatorname{dir}}

\newcommand{\cone}{\operatorname{cone}}

\newcommand{\rb}{\operatorname{rb}}
\newcommand{\myint}{\operatorname{int}}

\newcommand{\lspace}{\operatorname{lspace}}

\newcommand{\beq}{\begin{equation}}
\newcommand{\eeq}{\end{equation}}
\newcommand{\beqa}{\begin{eqnarray}}
\newcommand{\eeqa}{\end{eqnarray}}
\newcommand{\ba}{\begin{array}}
\newcommand{\ena}{\end{array}}
\newcommand{\bac}{\begin{array}{ccccccccccc}}
\newcommand{\eac}{\end{array}}
\newcommand{\bprop}{\begin{Proposition}}
\newcommand{\eprop}{\end{Proposition}}
\newcommand{\beqast}{\begin{eqnarray*}}
\newcommand{\eeqast}{\end{eqnarray*}}
\newcommand{\benum}{\begin{enumerate}}
\newcommand{\eenum}{\end{enumerate}}
\newcommand{\bit}{\begin{itemize}}
\newcommand{\eit}{\end{itemize}}
\newcommand{\bth}{\begin{Theorem}}
\newcommand{\enth}{\end{Theorem}}
\newcommand{\ble}{\begin{Lemma}}
\newcommand{\ele}{\end{Lemma}}
\newcommand{\bex}{\begin{Example}}
\newcommand{\eex}{\end{Example}}
\newcommand{\bcor}{\begin{Corollary}}
\newcommand{\ecor}{\end{Corollary}}
\newcommand{\brem}{\begin{Remark}}
\newcommand{\erem}{\end{Remark}}
\newcommand{\bass}{\begin{Assumption}}
\newcommand{\eass}{\end{Assumption}}
\newcommand{\bsmx}{\begin{small} \begin{pmatrix}}
\newcommand{\esmx}{\end{pmatrix} \end{small}}
\newcommand{\bpx}{\begin{pmatrix}}
\newcommand{\epx}{\end{pmatrix}}
\newcommand{\bbx}{\begin{bmatrix}}
\newcommand{\ebx}{\end{bmatrix}}
\newcommand{\bdef}{\begin{Definition}} 
\newcommand{\bconj}{\begin{Conjecture}}
\newcommand{\econj}{\end{Conjecture}}

\newcommand{\eps}{\epsilon} 
\newcommand{\commentout}[1]{}
\newcommand{\co}[1]{}

\newcommand{\nin}{\noindent}
\newcommand{\ti}{\times}
 
\newcommand{\pf}[1]{\vspace{.35cm} \nin {\bf Proof {#1} }} 
\newcommand{\norm}[1]{\parallel \! #1 \! \parallel}

\newcommand{\sym}[1]{{\cal S}^{#1}}
\newcommand{\psd}[1]{{\cal S}_+^{#1}}
\newcommand{\rad}[1]{\mathbb{R}^{#1}}

\newcommand{\eref}[1]{(\ref{#1})} 

\newcommand{\R}{ {\cal R} }

\newcommand{\la}{\langle}
\newcommand{\ra}{\rangle}

\newcommand{\tleft}{\triangleleft}
\newcommand{\tlefteq}{\trianglelefteq}

\newcommand{\LRA}{\Leftrightarrow}

\newcommand\emp{\emptyset}

\newcommand\xfh{x_{F,H}}
\newcommand\xfhi{x_{F,H_i}}

\title{\Large On the connection of facially exposed and nice cones}
\author{G\'{a}bor Pataki\thanks{Department of Statistics and Operations Research, University of North Carolina at Chapel Hill} \hspace{1cm} 
}
\date{}

\begin{document}

\maketitle 



\begin{abstract}  A closed convex cone $K$ in a finite dimensional Euclidean space 
is called {\em nice, } 
if the set $K^* + F^\perp$ is closed for all $F$ faces of $K,$ where 
$K^*$ is the dual cone of $K, \,$ and $F^\perp$ is the orthogonal complement of the linear span of $F.$ 
The niceness property plays a role in the facial reduction algorithm of Borwein and Wolkowicz,
and the question ``when is the linear image of the dual of 
a nice cone closed?''  also has a simple answer.

We prove several characterizations of nice cones and show a strong connection with 
facial exposedness. We prove that a nice cone must be facially exposed; in reverse, 
facial exposedness with an added condition implies niceness. 

We conjecture that nice, and facially exposed cones are actually the same, and give supporting 
evidence. 

\vspace{.2cm}

MSC2000 subject classification: Primary, 52A20, Secondary, 90C46, 49N15

\end{abstract} 




\maketitle

\pagestyle{myheadings}
\thispagestyle{plain}
\markboth{}{On the connection of facially exposed, and nice cones}

\section{Introduction}

Closed convex cones in finite dimensional Euclidean spaces 
appear in many areas of optimization. Conic linear programs -- optimization problems with 
a linear objective function and a feasible set expressed as the intersection of a 
closed convex cone with an affine subspace -- were introduced by Duffin in \cite{Duff:56}. 
They serve as a natural framework for studying
the duality theory of convex programs. The seminal interior-point framework 
of Nesterov and Nemirovskii \cite{NestNemirov:94} was also developed for conic LPs. 

The properties of the underlying cone determine, to a large extent, whether a conic LP is easy or hard. 
The nonnegative orthant is arguably the simplest cone useful in optimization. 
Second order, $p$-order, and semidefinite cones 
are more complex, but still admit efficient 
optimization algorithms 
(see e.g. \cite{AliGold:03}, \cite{Todd:00}, \cite{Guler:10}), 
and their geometry is also well understood (\cite{BarCar:75} and \cite[Appendix A]{Pataki:00A}). 
Copositive, and completely positive cones 
lie at the other end of the spectrum. Though they are very useful in optimization 
(see e.g. \cite{Burer:09, Duer:10}),
optimizing over them is more difficult. Also, while considerable progress has been made 
in describing their geometry (see \cite{AnstBurDo:09}, \cite{Dick:11}), a complete 
understanding (such as a complete description of their facial structure) 
is probably out of reach.

The goal of this paper is to study niceness, an intriguing geometric 
property of closed convex cones, and to connect it to facial exposedness. 
The niceness property is important for several reasons: first, 
it plays a role in the facial reduction algorithm of Borwein and Wolkowicz \cite{BorWolk:81}. 
Precisely, given a conic system $\{ \, x \, | \, g(x) \in K \, \}$ with $K$ a closed convex cone, 
the Borwein-Wolkowicz algorithm constructs a sequence of 
equivalent systems, the final one being strictly feasible. 
If $K$ is nice, then the reducing certificates can be chosen to be simpler in the algorithm.
For other aspects of facial reduction algorithms, we refer to \cite{WakiMura:10} and
\cite{PolikTerlaky:09}. 

Second, consider the following classical question: ``Is the linear image of a closed, convex cone closed?''
The closedness has a very simple characterization, when the {\em dual} of the cone is nice. 
First, let us note that for a convex set $C \,$ the relative interior of $C$ is denoted by $\ri C, \,$ 
for 
$x \in C \,$ the set of feasible directions at $x$ in $C$ is defined as
$\dir(x, C) \, = \, \{ \, y \, | \, x + \eps y \in C \, \text{for some} \, \eps > 0 \, \},$ 
and $\cl \dir(x,C)$ stands for the closure of $\dir(x,C).$ 
Also, for a linear map $M$ we denote by $\R(M)$ its rangespace, and by $M^*$ its adjoint map.
For motivation we recall a simplified version of Theorem 1.1 in \cite{Pataki:07}:
\bth \label{pataki-cl}
Let $M$ be a linear map, $C$ a nice cone, $C^*$ its dual cone, and $x \in \ri ( C \cap \R(M) ).$
Then 
\bit 
\item \label{mc} the set $M^*C^* \,$ is closed $\LRA$ $\R(M) \cap (\cl \dir(x, C) \setminus \dir(x,C)) = \emp.$ \qed
\eit 
\enth

(In Theorem \ref{pataki-cl} 
the set $\dir(x,C)$ is the same for all $x \in \ri ( C \cap \R(M) ), $ so the condition for the closedness of 
$M^* C^*$ only depends on $C$ and $M.$)

For better intuition, we can note that if $C$ is a polyhedral cone, then 
$\dir(x,C)$ is closed for all $x \in C,$ 
and that polyhedral cones are nice. So the direction $\Leftarrow$ above shows 
that $M^*C^*$ is closed for an arbitrary $M$ map, as expected. Also, if $x$ is in $\ri C \,$ 
(i.e., a ``Slater type'' condition is satisfied), 
then $\dir(x,C)$ is just a subspace, hence closed, so the same argument 
proves the closedness of $M^*C^*$ in this case as well. Thus Theorem \ref{pataki-cl} 
unifies two seemingly unrelated, sufficient conditions for the closedness of $M^*C^*.$ 

More recently, in \cite{GouveiaParriloThomas:12} Gouveia, Parrilo, and Thomas used 
the concept of niceness in studying the question 
whether a convex set can be represented as the projection of the intersection of a closed convex cone, and 
of an affine subspace. When the cone in question is nice, a sufficient condition for such a {\em lift} 
to exist becomes necessary and sufficient.

Facial exposedness of  convex cones is another classical concept in convex analysis. 
Many cones appearing in the optimization literature, for instance polyhedral, second order, 
$p$-cones, and the semidefinite cone  are both facially exposed, and nice: see for instance \cite{Pataki:07}. 

Here we study nice cones from two viewpoints: 
we describe characterizations (more precisely, we describe characterizations of the 
situation when $K^* + F^\perp$ is closed for a specific $F$ face of $K$), 
and find a direct, close connection with facial exposedness. 
In particular, we prove that a nice cone must be facially exposed; conversely,
facial exposedness with an added condition implies niceness. 
This leads us to raise the conjecture that the two classes of cones are actually the same, and 
to provide more supporting evidence.

The rest of the paper is structured as follows. In Section \ref{sect-prelim} we collect definitions 
and preliminary results. 
Section \ref{sect-char} has our main characterizations of nice cones, and describes the connection with facial exposedness. Section \ref{sect-conj} states the conjecture, shows a supporting example, 
and shows that proving a seemingly weaker version would already suffice. 
In this section we also describe another characterization of nice cones, and shows how it 
may lead to the proof of the main conjecture.

\section{Preliminaries} \label{sect-prelim} 

Throughout the paper we assume that the underlying space is a finite dimensional Euclidean space. 
For a set $S$ we write $\cl S$ for its closure, $\lin S$ for its linear span,
and $S^\perp$ for the orthogonal complement of its linear span. 
For a convex set $S$ we denote its relative interior by $\ri S,$ and its relative boundary by $\rb S.$ 
For a one-element set $\{ y \}$ we abbreviate $\{ y \}^\perp$ by $y^\perp.$  

A set $C$ is called a  {\em cone}, if $\lambda x \in C$ holds for all $x \in C, \,$ and $\lambda \geq 0.$ 
For a set $S$ the set of all nonnegative combinations of elements of 
$S$ is clearly a convex cone, which 
is called the {\em cone generated by} $S,$ and denoted by 
$\cone S.$ For a one-element set $\{ y \}$ we abbreviate $\cone \{ y \}$ by $\cone y. \,$ 

General references on convex analysis that we used are
for instance \cite{Rockafellar:70, BorLewis:00, Guler:10}. References \cite{Berman:73, Tam:85, Barker:73} 
cover more specifically the theory of cones. 
If $C$ is a convex cone in a Euclidean space $X, \,$ then its {\em lineality space} is defined as 
$$
\lspace C \, = \, C \cap - C, 
$$
and its dual cone as 
$$
C^* \, = \, \{ \, y \in X \, | \, \la y, x \ra \geq 0 \,\, \forall x \in C \}.
$$
We say that $C$ is {\em pointed}, if $\lspace C = \{0 \}.$ 
For convex cones $C, C_1, \,$ and $C_2$ we have 
\beqa
C^{**} & = & \cl C, \\
(C_1 + C_2)^* & = & C_1^* \cap C_2^*.
\eeqa
Furthermore, if $C_1$ and $C_2$ are also closed, then 
\beqa \label{k1k2*}
(C_1 \cap C_2)^* & = & \cl (C_1^* + C_2^*).
\eeqa
Given a closed convex cone $C, \,$ and $x_1, x_2 \in C, \,$ the open line-segment between $x_1$ and $x_2$ is defined 
as 
$$
]x_1, x_2[ \, = \, \{ \, \lambda x_1 + (1-\lambda) x_2 \, | \, 0 < \lambda < 1 \, \}.
$$
A convex subset $E$ of $C$ is called a {\em face of } $C, \,$ if 
$x_1, x_2 \in C, \, ]x_1, x_2[ \cap E \neq \emptyset$ implies that $x_1$ and $x_2$ are both in $E;$ 
equivalently, if $x_1 + x_2 \in E$ imply that $x_1$ and $x_2$ both are in $E.$ 
The cone $C$ itself is clearly a face of $C, $ and all faces of $C$ are cones in their own right. 
For all $x \in C $ there is a unique minimal face of $C$ that contains $x, \,$ namely the 
face having $x$ in its relative interior. 

We write $E \tlefteq C$ to denote that $E$ is a face of $C, \,$ and $E \tleft C$ to abbreviate 
\mbox{$E \tlefteq C, \, E \neq C.$}
The definition implies that the intersection of faces is again a face. Also, $\lspace C$ is the inclusionwise minimal face of $C.$ 
It is straightforward to show that if $E_1 \tlefteq C, \,$ and $E_2 \subseteq E_1, \,$ then 
$E_2$ is a face of $E_1$ iff it is a face of $C.$ 

We call a face $E_1$ of $C$ a
{\em properly maximal face of} $C, \,$ if $E_1 \neq C, \,$ and there is no $E_2$ such that 
$E_1 \tleft E_2 \tleft C.$ 
We call an $E_1$ face of $C$ a {\em properly minimal face of} $C,\,$ if $E_1 \neq \lspace C, \,$ and 
there is no $E_2$ such that $\lspace C \tleft E_2 \tleft E_1.$ 
Properly minimal faces of a pointed, closed convex
cone $C$ are of the form $\{ \, \lambda x \, | \, \lambda \geq 0 \, \}, \,$ 
where $x \in C \setminus \{ 0 \}, \,$ and are called {\em extreme rays.} 

For example, if $C$ is the the nonnegative orthant in $\rad{n}, \,$ then its properly maximal faces 
are $E_i = \{ \, x \, | \, x \in C, \, x_i = 0 \, \}$ for $i=1, \dots, n, \,$ and its properly minimal 
faces are cones generated by unit vectors. 
If $C$ is a halfspace, 
i.e., $C = \{ \, x \, | \, \la a, x \ra \geq 0 \,\} \,$ for some $a \neq 0, \,$ then its only properly
maximal face is its lineality space $\{ \, x \, | \, \la a, x \ra = 0 \, \}, $ and 
its only properly minimal face is $C$ itself.

A remark on notation: we will look at characterizations of the niceness of a closed convex cone, and 
will generally denote this cone by $K.$ In collecting relevant results we usually 
reference  a closed convex cone by $C, \,$ since the role of $C$ later on will be played sometimes by $K, \,$ and sometimes by $F^*, \,$ where $F$ is a face of $K.$

A subset $E$ of $C$ is called an {\em exposed face of} $C, \,$ if it is the intersection of $C$ with a 
supporting hyperplane, i.e.,
$$
E \, = \, C \cap y^\perp
$$
for some $y$ satisfying $\la y, x \ra \geq 0$ for all $x \in C, \,$ i.e., $y$ must be in $C^*.$ 
We say that $y$  {\em exposes} $E$. 
Also, if $H$ is the smallest face of $C^*$ that contains $y, \,$ then $E = C \cap H^\perp$ 
holds for the above $E.$ 

An exposed face of $C$ is always a face, but a face $E_1$ may not be exposed. This happens when every 
supporting hyperplane of $C$ that contains $E_1$ actually contains a larger face $E_2, \,$ 
i.e., there is an $E_2$ face of $C$ with $E_1 \tleft E_2, \,$ such that 
$E_1 \subseteq y^\perp$ implies $E_2 \subseteq y^\perp$ for all $y \in C^*.$ 
An equivalent statement is that 
\co{Since for $y \in C^*$ and for $i=1,2, \,$ we have $E_i \subseteq y^\perp$ iff 
$y \in E_i^\perp, \,$ an equivalent statement is }
$ C^* \cap E_2^\perp = C^* \cap E_1^\perp \,$ holds (with the containment $\subseteq$ being trivial).
Example \ref{ex1} shows a cone with a nonexposed face.

We say that a closed convex cone $C$ is {\em facially exposed}, if all of its faces are exposed. Based on the above 
argument, an equivalent definition is requiring 
\beq
C^* \cap E_2^\perp \subsetneq C^* \cap E_1^\perp \; 
\eeq
for all $E_1$ and $E_2$ faces of $C$ with $E_1 \subsetneq E_2.$ 

The intersection of exposed faces is again an exposed face, so the smallest exposed face containing a 
subset of $C$ is well-defined. 
In particular, if $E$ is a face of $C, \,$ then the smallest exposed face containing it is 
$C \cap H^\perp \, = \, C \cap y^\perp, \,$ where $H = C^* \cap E^\perp, \,$  and $y \in \ri H. \,$ 
\co{or, equivalently, to 
$C \cap H^\perp, $ where $H = C^* \cap E^\perp.$ }
Hence a face $E$ is exposed, iff it is equal to the smallest exposed face containing it. 

The following proposition is well-known -- see for instance 
Proposition 2.1 in \cite{Tam:85}. 
\bprop \label{properprop}
Suppose $C$ is a closed convex cone, and $E \tlefteq C.$ Then 
\benum
\item $E = C             \, \Leftrightarrow \, C^* \cap E^\perp = C^\perp.$ 
\item $E = \lspace C \, \Leftrightarrow \, C^* \cap E^\perp = C^*.$
\eenum
\eprop
\qed

The space of $n$ by $n$ symmetric, 
and the cone of $n$ by $n$ symmetric, positive semidefinite  matrices are  denoted by $\sym{n}$, 
and $\psd{n}$, respectively. 
The space $\sym{n}$ is equipped with the inner  product 
\beqast
 X \bullet Z  & := & \sum_{i,j=1}^n x_{ij} z_{ij}, 
\eeqast
where the components of $X$ and $Z$ are denoted by $x_{ij}$ and $z_{ij}, \,$ respectively, 
and it is a well-known fact, that $\psd{n}$ is self-dual with respect to this inner product.

The faces of $\psd{n}$ have an attractive, and simple description.
After applying a rotation $V^T (.) V$ by a full-rank matrix $V, \,$ 
any face can be brought to the form 
\beqast
E & = & \Biggl\{ \bpx X & 0 \\ 0 & 0 \epx \, | \,  X \in \psd{r} \, \Biggr\}.
\eeqast
For a proof, see \cite{BarCar:75}, or Appendix A in \cite{Pataki:00A} for a somewhat simpler one.
(Exposed faces of more general spectral sets, with the  semidefinite cone being 
a special case, have been characterized in \cite{Lewis:98}.) 
For a face of this form we will use the shorthand
\beq \label{f-oplus}
E = \begin{pmatrix}  \oplus \! &  0  \\
                    0   \! &  0 
     \end{pmatrix}, \, 
\lin E = \begin{pmatrix}  \ti \! &  0  \\
                    0   \! &  0 
     \end{pmatrix}, \, 
\eeq
when the size of the partition is clear from the context. 
The $\oplus$ sign denotes a positive semidefinite submatrix, and the sign
 $\ti$ stands for a submatrix with arbitrary
elements. We use similar notation for other subsets of $\sym{n}:$ for instance, 
$$
\begin{pmatrix}  \oplus \! &  \ti   \\
                   \ti    \! &  \ti  
     \end{pmatrix}
$$
stands for the set of matrices with the upper left block positive semidefinite, and the other elements 
arbitrary. 

Facially nonexposed cones can be constructed by taking sums of facially exposed ones, as 
Example \ref{ex1} shows. The cross-section of the cone in Example \ref{ex1} is
illustrated on Figure \ref{fig1}. We give this example in detail, since we will return to it later.
\compilepdf{
\begin{figure}[h]
\begin{center}
\includegraphics[width=4in]{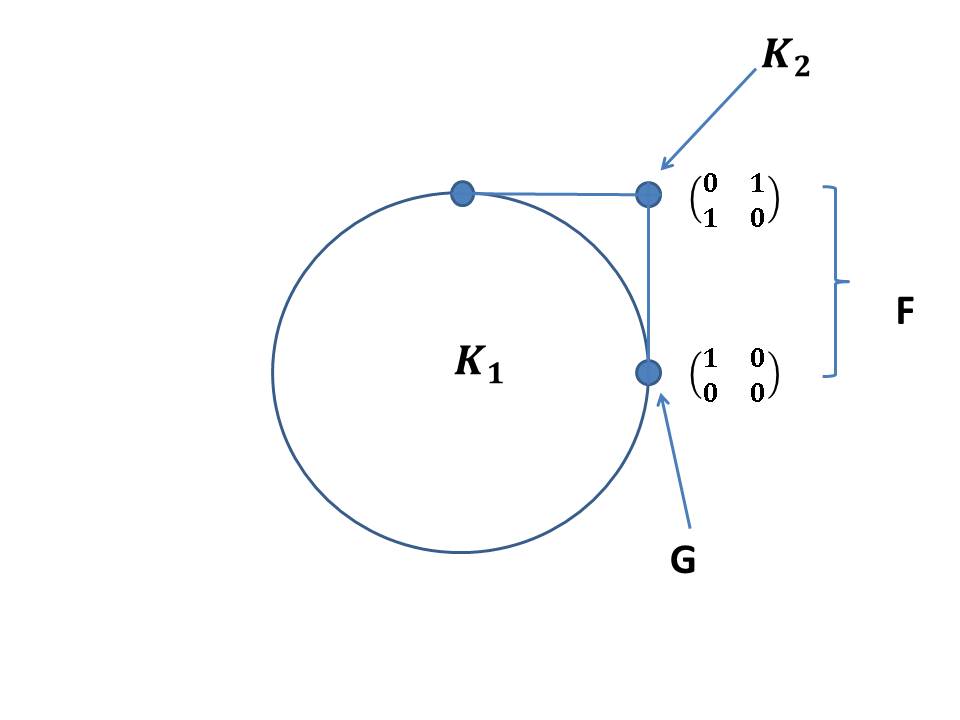}   
\caption{\small Cross section of a facially nonexposed cone \normalsize}
\label{fig1}
\end{center}
\end{figure}
}
\bex \label{ex1}
{\rm 

Define the cone $K \subseteq \sym{2}$ as $K = K_1 + K_2, \,$ with $K_1 = \psd{2}, \, K_2 = \cone \bpx 0 & 1 \\ 1 & 0 \epx.$ 
Let 
\beqa
G & = & \cone \bpx 1 & 0 \\ 0 & 0 \epx, \\
F & = & \cone \{ G \cup K_2 \}. 
\eeqa
It is straightforward to check that $K$ is closed, and that $G$ and $F$ are both faces of $K.$ 

Also, 
\beq \label{k*ex1-1} 
K^* \, = \, K_1^* \cap K_2^* \, = \, \{ \, X \in \psd{2} \, | \, \, x_{12} \geq 0 \, \}.
\eeq
Hence $K^* \cap G^\perp \, = \, K^* \cap F^\perp \, = \, \cone Y, $
where 
$$
Y \, = \, \bpx 0 & 0 \\ 0 & 1 \epx,
$$
so $G$ is not exposed. Clearly, $Y \in \ri (K^* \cap G^\perp), \,$ 
so the smallest exposed face of $K$ that contains $G$ is 
$K \cap Y^\perp = F, \,$ and this fact can also easily be checked 
by looking at Figure \ref{fig1}.
}
\eex

We repeat the  main definition of the paper for convenience: 
\bdef \label{nice-def}
A closed convex cone $K$ is called nice, if the set $K^* + F^\perp$ is closed for all $F \tlefteq K.$
\end{Definition}
\brem \label{nice-def-remark}
{\rm 
Since for a closed convex cone $K, \,$ and $F \tlefteq K$ we have 
$F = K \cap \lin F,$ by \eref{k1k2*} it follows that $F^* = \cl (K^* + F^\perp). \,$ 
Hence Definition \ref{nice-def} is equivalent to requiring 
\beq \label{FK}
F^* \, = \, K^* + F^\perp \,\, {\rm for \, all} \,\, F \tlefteq K.\
\eeq
Also, \eref{FK} trivially holds for $F = \lspace K, \,$ and $F = K, \,$ so it suffices to require it for
other  faces of $K$ in Definition \ref{nice-def}. 
}
\erem
Parts of the following proposition, which collects properties 
of closed, convex, possibly nonpointed cones are 
based on the remarks in Section 18 in \cite{Rockafellar:70}, and the rest are easy exercises to prove.
\bprop \label{lspace-prop}
Let $C$ be a closed convex cone, and $L = \lspace C.$ Then the following statements hold:
\benum
\item $C = C_0 + L, \,$ where $C_0 = C \cap L^\perp,$ and $C_0$ is pointed.
\item \label{part-later} If $E \tlefteq C, \,$ then $\lspace E  =  \lspace C.$
\item The mapping $E \rightarrow E_0 := E \cap L^\perp, $ where $E \tlefteq C, \,$ 
defines a one-to-one correspondence 
between the faces of $C$ and $C_0;$ in particular, $E = E_0 + L. \,$ 
\item \label{EE0} For the above $E$ and $E_0$ it holds that 
\benum 
\item \label{EE00} $E$ is a properly minimal face of $C$ iff $E_0$ is an extreme ray of $C_0.$ 
\item \label{again} $E$ is exposed iff $E_0$ is.
\eenum
\eenum
\qed
\eprop
In general it seems to be known that a nonpointed, closed, convex cone is generated by the union of its 
properly minimal faces. Since we were not able to find a result stated precisely in this form, we 
state, and prove:
\bprop \label{CE}
Let $C$ be a closed convex cone.
Then 
\beq \label{f*h}
C \, = \, \cone \, \bigcup \, \{ \, E \, | \, E \, {\rm is \, a \, properly \, minimal \, face \, of \,} C \, \}.
\eeq
\eprop
\pf{} Let $L = \lspace C, \,$ and write $C = C_0 + L, \,$ with $C_0 = C \cap L^\perp.$ 
Then 
$$
\begin{array}{rcl}
C & = & \cone \, \bigcup \, \{ \, E_0 \, | \, E_0 \, {\rm is \, an \, extreme \, ray \, of \,} C_0 \, \} + L \\
    & = & \cone \, \bigcup \, \{ \, E \cap L^\perp \, | \, E \, {\rm is \, a \, properly \, minimal \, face \, of \,} C \, \} + L \\
    & = & \cone \, \bigcup \, \{ \, E \cap L^\perp + L \, | \, E \, {\rm is \, a \, properly \, minimal \, face \, of \,} C \, \}  \\
    & = & \cone \, \bigcup \, \{ \, E \, | \, E \, {\rm is \, a \, properly \, minimal \, face \, of \,} C \, \}.
\end{array}
$$
Here the first equation 
comes from the fact that $C_0$ is pointed, and Theorem 18.5 in \cite{Rockafellar:70}, 
the third is trivial, and the others follow from Proposition \ref{lspace-prop}. 
\qed

\section{Characterizations of nice cones, and connections to facial exposedness}
\label{sect-char}

Throughout this section we assume that 
\begin{center}
\framebox{ $K$ is a closed, convex cone.}
\end{center}

In Theorem \ref{nice-thm1} and Remark \ref{rem-nice1} 
we give several characterizations of the situation when $F^* = K^* + F^\perp$ 
holds for a specific $F$ face of $K.$ In Theorem \ref{nice-thm2} we build on this to make the connection of
the niceness of $K$ to its facial exposedness.

We start with an informal discussion. 
If $F$ is a face of $K, \,$ then $\lspace F^* = F^\perp$ holds, hence a face $H$ of $F^*$ satisfies 
$ H \supseteq F^\perp.$ 
Clearly,
\beq \label{k*h}
K^* \cap H \supseteq K^* \cap F^\perp,
\eeq
and by the definition of faces, and $K^* \subseteq F^*$ both sets in \eref{k*h} are faces of $K^*.$ 
However, they may be equal, even when $H$ and $F^\perp$ are not. 

The equality of $F^*$ and $K^* + F^\perp$ is characterized by strict containment holding in \eref{k*h} for 
{\em all} $H$ properly minimal faces of $F^*$ (i.e., minimal faces that are distinct from $F^\perp$);
equivalently, by strict containment holding for {\em all} $H$ faces that are distinct from 
$F^\perp.$ These are conditions \eref{threep} 
in Theorem \ref{nice-thm1}, and (3') in Remark \ref{rem-nice1}, and we use them later 
in Theorem \ref{nice-thm2} to connect niceness to facial exposedness.

We need the following result: 
\bprop \label{Hprop}
Let $F \tlefteq K,$ and $H \tlefteq F^*.$ 
Then 
\beqa \label{lspaceH}
\lspace H & = & F^\perp, \\ \label{rilspaceH}
\ri H     & = & \ri H + F^\perp. 
\eeqa
\eprop
\pf{} Statement \eref{lspaceH} directly follows from 
part \eref{part-later} in Proposition \ref{lspace-prop}. Statement 
\eref{rilspaceH} comes from 
\eref{lspaceH} and 
the fact that $H$ is a closed convex cone in its own right, hence it is easy to check 
that $\ri H = \ri H + \lspace H.$ 
\qed

\bth
\label{nice-thm1} 
Let $F \tlefteq K.$ Then the following statements are equivalent:
\benum
\item \label{one} $F^* = K^* + F^\perp.$ 
\item \label{twop} $K^* \cap \ri H \neq \emptyset$ holds for all $H \tlefteq F^*.$
\item \label{threep} $K^* \cap H \supsetneq K^* \cap F^\perp$ holds for all $H$ properly minimal faces of $F^*.$ 
\eenum
\enth
\pf{of $\eref{one} \Rightarrow \eref{twop}$} Let $H$ be a face of $F^*, \,$ and $x \in \ri H.$ Write 
$x = x_1 + x_2, \,$ with $x_1 \in K^*, \, x_2 \in F^\perp.$ 
Hence $x_1 = x - x_2 \in \ri H + F^\perp = \ri H, \,$ where the last equation follows from 
\eref{rilspaceH}. 

\pf{of $\eref{twop} \Rightarrow \eref{threep}$} This implication 
follows from the fact that if $H$ is a properly minimal face of $F^*, \,$ then its only face other than itself is
$F^\perp, \,$ hence $\rb H = F^\perp.$ 

\pf{of $\eref{threep} \Rightarrow \eref{one}$} Proposition \ref{CE} implies 
\beq \label{f*h2}
F^* \, = \, \cone \, \bigcup \, \{ \, H \, | \, H \, {\rm is \, a \, properly \, minimal \, face \, of \,} F^* \, \}.
\eeq
Let $H$ be an arbitrary properly minimal face of $F^*, \,$ and assume that 
$K^* \cap H \supsetneq K^* \cap F^\perp \,$ holds. 
Given \eref{f*h2}, it suffices to prove $H \subseteq K^* + F^\perp.$ 
As remarked above, $H$ is the disjoint union of $\ri H$ and $F^\perp.$

Let $x \in H.$ If $x \in F^\perp, \,$ then of course $x \in K^* +  F^\perp, \,$ so suppose 
$ x \not \in F^\perp,$ i.e., $x \in \ri H.$ 
By the assumption there is $y \in (K^* \cap H ) \setminus F^\perp. \,$ 
If $ x = y, \,$ then again $x \in K^* + F^\perp.$ If $x \neq y, \,$ then let us define the two half-lines 
$$
\begin{array}{rcl}
r^+_{x,y} & = & \{ \, x + \lambda y \, | \, \lambda \geq 0 \, \}, \\
r^-_{x,y} & = & \{ \, x - \lambda y \, | \, \lambda \geq 0 \, \}.
\end{array}
$$
Then $r^+_{x,y} \subseteq H, \,$ since $H$ is a convex cone, and $x$ and $y$ are both in $H.$ 
Hence $r^-_{x,y} \not \subseteq H, \,$ since both 
$r^+_{x,y}$ and $r^-_{x,y}$ being in $H$ would imply $ y \in \lspace H = F^\perp.$ 
Define 
$$
\lambda^* \, = \, \max \, \{ \, \lambda \, | \, x - \lambda y \in H \, \}.
$$
Since $r^-_{x,y} \not \subseteq H, \,$ we have $\lambda^* < + \infty;$ 
since $H$ is closed, $\lambda^*$ is attained, and by $x \in \ri H$ we have $\lambda^* > 0.$ 
Let $z = x - \lambda^* y.$ Clearly, $z$ must be in the relative boundary of $H, \,$ i.e., 
$z \in F^\perp.$ 
Hence $x = \lambda^* y + z \, \in \, K^* + F^\perp, $ as required.

\qed

\brem \label{rem-nice1}
{\rm 

If $F$ is as in Theorem \ref{nice-thm1}, it is straightforward to see that two other conditions equivalent to 
$F^* = K^* + F^\perp$ are 
\benum
\item[(2')] \label{three} $K^* \cap \ri H \neq \emptyset$ holds for all $H$ properly minimal faces of $F^*.$
\item[(3')] \label{two} $K^* \cap H \supsetneq K^* \cap F^\perp$ holds for all $H \tlefteq F^*$ s.t. 
$H \neq F^\perp.$ 
\eenum
Indeed, it is easy to check that (using the numbering of statements in Theorem \ref{nice-thm1}), that the 
implications 
$(2) \Rightarrow (2') \Leftrightarrow (3) \,$ and $(2) \Rightarrow (3') \Rightarrow (3) \,$ hold. 

Also, for an $H$ face of $F^*,$ we have $\ri H = \ri H + F^\perp$ by \eref{rilspaceH}. 
Hence 
\beq \label{k*hequiv}
K^* \cap \ri H \neq \emptyset \LRA K^* \cap (\ri H + F^\perp) \neq \emptyset \LRA (K^* + F^\perp) \cap \ri H \neq \emptyset,
\eeq
so replacing $K^*$ by $K^* + F^\perp$ in \eref{twop} in Theorem \ref{nice-thm1} and 
(2') above yields equivalent conditions. 

Also, since $\lspace H = F^\perp, \,$ it is easy to check that 
$$ 
K^* \cap H \supsetneq K^* \cap F^\perp \LRA (K^* + F^\perp) \cap H \supsetneq (K^* + F^\perp) \cap F^\perp
$$ 
(and the last set is just $F^\perp$). Thus, replacing 
$K^*$ by $K^* + F^\perp$ in 
\eref{threep} in Theorem \ref{nice-thm1} and (3') above we also obtain equivalent conditions.
}
\erem

\bth
\label{nice-thm2} 
The following statements hold.
\benum
\item \label{nice2-1} If $K$ is nice, then it is  facially exposed.
\item \label{nice2-3} If $K$ is facially exposed, and for all $F \tlefteq K$  all 
properly minimal faces of  $F^*$ are exposed, then $K$ is nice.
\eenum
\enth

\pf{} Consider the statements
\beq \label{a}
K^* \cap H \supsetneq K^* \cap F^\perp, \,
\eeq
where $F \tlefteq K \,$ and $H$ is a face of $F^*$ distinct from $F^\perp, \,$ 
and 
\beq \label{b}
K^* \cap G^\perp \supsetneq K^* \cap F^\perp, 
\eeq
where $F$ and $G$ are faces of $K$ satisfying $G \subsetneq F.$ 

Theorem \ref{nice-thm1} and Remark \ref{rem-nice1} show that $K$ is nice, iff \eref{a} holds 
for {\em all} 
$F \tlefteq K \,$ and {\em all} $H$ properly minimal faces of $F^*;$ or equivalently, 
for all $F \tlefteq K, \,$ and {\em all} $H$ faces of $F^*$ that are distinct from $F^\perp. \,$ 
Also, $K$ is facially exposed, iff \eref{b} holds for {\em all} $F$ and $G$ faces of $K$ with $G \subsetneq F.\,$ 

\co{To prove \eref{nice2-1}, let $F$ and $G$ be faces of $K$ with $G \subsetneq F.$ 
Let us define $H = F^* \cap G^\perp.$ 
Since $G \neq F,$ Proposition \ref{properprop} implies $H \neq F^\perp.$ 
As $K$ is nice, \eref{a} holds, and since
\beq \label{k*h2}
K^* \cap H = K^* \cap F^* \cap G^\perp = K^* \cap G^\perp,
\eeq
\eref{b} follows. 
}
To prove \eref{nice2-1}, assume that $K$ is nice, and let $F$ and $G$ be faces of $K$ with $G \subsetneq F.$ 
We will prove that \eref{b} holds.
Let us define $H = F^* \cap G^\perp.$ 
Since $G \neq F,$ Proposition \ref{properprop} implies $H \neq F^\perp.$ 
As $K$ is nice, \eref{a} holds, and since
\beq \label{k*h2}
K^* \cap H = K^* \cap F^* \cap G^\perp = K^* \cap G^\perp,
\eeq
\eref{b} follows. 

To prove \eref{nice2-3}, assume that the condition therein is satisfied, 
let $F$ be a face of $K, \,$ and $H$ a properly minimal face of $F^*.$ 
We will prove that \eref{a} holds.
By the assumption $H$ is an exposed face, so $H = F^* \cap G^\perp$ holds for a $G$ face of $F.$ Since 
$H \neq F^\perp, \,$ by Proposition \ref{properprop} we have $G \neq F.$ 
Then clearly \eref{k*h2} holds. Since $K$ is facially exposed, \eref{b} holds as well, hence \eref{a} follows.

\qed

To better understand cones that are {\em not} nice, 
we will look at $F$ faces of $K \,$ s.t. $K^* + F^\perp$ is not closed, when there is such a face, 
i.e., (cf. Remark \ref{nice-def-remark}), when the set $F^* \setminus (K^* + F^\perp)$ is nonempty.
The following corollary shows how to find points in this difference set.
\bcor \label{corK} The following statements hold.
\benum
\co{\item \label{1} If $F \tlefteq K, \,$ and $H$ is a face of $F^*$ s.t. 
$H \neq F^\perp, \,$ and 
$K^* \cap H = K^* \cap F^\perp, $ then 
$$
\ri H \subseteq F^* \setminus (K^* + F^\perp).
$$}
\item \label{1} If $F \tlefteq K, \,$ $H$ is a face of $F^*$ distinct from $F^\perp, \,$ and 
$K^* \cap H = K^* \cap F^\perp, $ then 
$$
\ri H \subseteq F^* \setminus (K^* + F^\perp).
$$
\item \label{2} If $G$ is a nonexposed face of $K, \,$ and $F$ the smallest exposed face of $K$ that contains $G, \,$ then 
$$
\ri (F^* \cap G^\perp) \subseteq F^* \setminus (K^* + F^\perp).
$$
\eenum
\ecor
\pf{of \eref{1}} The containment $\ri H \subseteq F^*$ is obvious. 
Since $F^\perp \tleft H,$ we have $\ri H \cap F^\perp = \emptyset, \,$ and this 
with $K^* \cap H = K^* \cap F^\perp $ implies 
$K^* \cap \ri H = \emptyset. \,$ 
In turn, the equivalence \eref{k*hequiv} proves 
$(K^* + F^\perp) \cap \ri H = \emptyset.$ 

\pf{of \eref{2}} Let us define $H = F^* \cap G^\perp.$ 
Since $G \neq F, \,$ by Proposition \ref{properprop} we obtain $H \neq F^\perp.$
Since $G$ is a nonexposed face of $K$, 
and $F$ is the smallest exposed face that contains it, 
we have $K^* \cap F^\perp = K^* \cap G^\perp, \,$ 
hence $K^* \cap H = K^* \cap F^\perp, \,$ so part \eref{1} implies our claim.

\qed

\noindent{\bf Example \ref{ex1} continued} With $G$ a nonexposed face of $K, \,$ and $F$ the smallest exposed
face containing it, we have 
$$
\ba{rcl}
F^*              & = &  \{ \, X \in \sym{2} \, | \, \, x_{11} \geq 0, \, x_{12} \geq 0 \, \}, \\
F^* \cap G^\perp & = &  \{ \, X \in \sym{2} \, | \, \, x_{11} = 0, \, x_{12} \geq 0 \, \}, \\
F^\perp          & = &  \{ \, X \in \sym{2} \, | \, \, x_{11} = 0, \, x_{12} = 0 \, \}.
\ena
$$
Let 
$$
X = \bpx 0 & 1 \\ 1 & 0 \epx.
$$
Then clearly $X \in \ri ( F^* \cap G^\perp ), \,$ hence Corollary 
\ref{corK} implies $X \not \in K^* + F^\perp.$ 
Given the description of $K^*$ in \eref{k*ex1-1} one can indeed easily verify this fact.

\section{Are facially exposed, and nice cones the same?}
\label{sect-conj}

The main conjecture of the paper is:

\bconj \label{conj1}
A closed convex cone is nice if and only if it is facially exposed. \qed
\econj

Proving Conjecture \ref{conj1} would be very interesting, since facial exposedness and niceness 
are both fundamental, and at first sight unrelated geometric properties of cones.

Theorem \ref{nice-thm2} already finds a strong connection: niceness implies facial exposedness,
and facial exposedness with an added condition implies niceness. In support of 
Conjecture \ref{conj1}, we first present an example to show that the added condition in general is not 
necessary. Next, in Theorem \ref{weaker-thm} 
we show that proving a weaker version of Conjecture \ref{conj1} would already be sufficient.
Finally, we give a different characterization of nice cones in 
Corollary \ref{cor-xfh}, and outline how 
this may lead to a proof of Conjecture \ref{conj1}.

We need the following result:
\bprop \label{chua} Suppose that $K_1$ and $K_2$ are nice cones. Then $K_1 \cap K_2$ is also nice.
\eprop
\pf{} We will use a result of Chua and Tun\c{c}el in \cite{ChuaTuncel:08}. 
First, for a set $S, \,$ and a map $L, \,$ let us define $L^{-1}(S)$ 
as the preimage of $S$ under the map $L, \,$ 
i.e., 
$
L^{-1}(S) \, := \, \{ \, x \, | \, L(x) \in S \, \}.
$
In \cite{ChuaTuncel:08} a closed convex cone $K,$ which is pointed, and has nonempty interior is called 
{\em $G$-representable}, if $G$ is a cone of the same type, and there exists a linear map $L$ such that
$$
\myint K = L^{-1} (\myint G).
$$
Theorem 6.7 in \cite{Rockafellar:70} implies that this is equivalent to 
$$
K = L^{-1} (G),
$$
and Proposition 4 in \cite{ChuaTuncel:08} shows that if $G$ is nice, then so is $K.$ 
In fact, it is not hard to slightly modify Proposition 4 in \cite{ChuaTuncel:08} to show that 
if $K = L^{-1}(G)$ for a linear map $L, \,$ 
and $G$ is nice, then so is $K, \,$ i.e., we do not have to assume 
pointedness and full-dimensionality of $G$ and $K.$ 

Now suppose that $K_1$ and $K_2$ are nice cones, define the cone $G$ as $G = K_1 \ti K_2, \,$ 
and the linear map $L$ as $L(x) = (x,x).$ Then clearly $G$ is nice, and 
$L^{-1}(G) = K_1 \cap K_2,$ hence our claim follows.
\qed

Precisely, Example \ref{ex2} shows a closed, convex, facially exposed cone $K, \,$ 
which is nice, however, there is a face $F$ of $K \,$ such that 
an $H$ properly minimal face of $F^* \,$  is not exposed. We first informally describe Example \ref{ex2}.
We construct $K$ as $K = K_1 \cap K_2, \,$ where $K_1$ is a semidefinite cone, and $K_2$ is a halfspace. 
By Proposition \ref{chua} we have that $K_1 \cap K_2$ is nice. The cones $K_1$ and $K_2$ are also chosen 
so that their relative interiors intersect, hence 
(see e.g. Section 5 in \cite{Pataki:07}) $K^* = K_1^* + K_2^*.$ 

Then we choose suitable faces $F_1$ of $K_1, \,$ and $F_2$ of $K_2.$ 
The definition of faces implies that $F := F_1 \cap F_2$ is a face of $K$ 
(in fact, a theorem of Dubins in \cite{Dubins:62} shows that all faces of $K$ arise 
in this manner). 
Also, $F_1$ and $F_2$ are chosen to satisfy $\ri F_1 \cap \ri F_2 \neq \emp, \,$ so $F^* = F_1^* + F_2^*.$ 
As $F^*$ is the sum of two simple, facially exposed closed, convex cones, 
one can expect it to have nonexposed faces, like $K$ does in Example \ref{ex1}, and 
we can rigorously show that there is indeed such a face, which is properly minimal.

\bex
\label{ex2}
{\rm 
Let 
$$
M \, = \, \bpx 0 & 1 & 0 \\ 1 & 0 & 0 \\ 0 & 0 & 0 \epx,
$$
and define $K = K_1 \cap K_2, \,$ with 
\beq
K_1 \, = \, \psd{3}, \, K_2 \, = \, \{ \, X \in \sym{3} \, | \, M \bullet X \geq 0 \, \}.
\eeq
Also define $F \, = \, F_1 \cap F_2,$ 
with 
$$
F_1 \, = \, 
\bpx 
  \parbox{1cm}{$\,\, \bigoplus$}   & \hspace{-.5cm} \parbox{1cm} {$\ba{c} 0 \\ 0 \ena$} \hspace{-.5cm} \\
  \parbox{1cm}{$0 \,\,\,\, 0$}  &   \hspace{-.4cm} 0 
\epx, \,\, F_2 = K_2.
$$
(Recall the notation for the faces of the semidefinite cone from Section \ref{sect-prelim}.)

Then $\ri K_1 \cap \ri K_2 \neq \emptyset, \,$ hence $K^* = K_1^* + K_2^*, \,$ where 
$$
K_1^* \, = \, \psd{3}, \, K_2^* = \cone M.
$$
It is easy to check that $\ri F_1 \cap \ri F_2 \neq \emp, \,$ hence $F^* = F_1^* + F_2^*, \,$ where 
$$
F_1^* \, = \, 
\bpx 
  \parbox{1cm}{$\,\, \bigoplus$}   & \hspace{-.5cm} \parbox{1cm} {$\ba{c} \ti \\ \ti \ena$} \hspace{-.5cm} \\
  \parbox{1cm}{$\ti \,\,\,\, \ti$}  &   \hspace{-.4cm} \ti 
\epx, \,\, F_2^* = K_2^* = \cone M.
$$
(More formally, $F_1^*$ is the set of $3$ by $3$ symmetric matrices, whose upper left $2$ by $2$ block is 
positive semidefinite, and the rest of the components are arbitrary.) 
Now, let us define
$$
H \, = \, \bpx 0 & 0       & \ti \\
               0 & \oplus  & \ti \\
             \ti & \ti     & \ti
            \epx.
$$
}
\eex
(Again, more formally $H$ is the set of $3$ by $3$ symmetric matrices 
$X$ with $x_{11} = x_{12} = x_{21} = 0, \,$ 
$x_{22} \geq 0, \,$ and the rest of the components arbitrary.) 

\bprop
If $K, F, \,$ and $H$ are as in Example \ref{ex2}, then 
$H$ is a properly minimal face of  $F^*,$ which is not exposed.
\eprop
\pf{} We first prove that $H$ is a face. Let $X \in H, \,$ and suppose $X = Y+Z, \,$ where $Y, Z \in F^*.$ 
We show that $Y$ and $Z$ are in $H.$ 

Since $F^* = F_1^* + F_2^*, \,$ 
we can write $Y = S+T, \, Z = U + V, \,$ where \mbox{$S, U \in F_1^*, \,$}
\mbox{$T, V \in F_2^*.$}
Let us write $s_{ij}, \, t_{ij}, \, u_{ij}, \, v_{ij}$ for the components of $S, T, U, \,$ and $V, \,$ respectively. Since 
$t_{11} = v_{11} = 0, \,$ we have 
\beq \label{x11}
x_{11} \, = \, s_{11} + u_{11}.
\eeq
With $x_{11} = 0, \, s_{11} \geq 0, \, u_{11} \geq 0, \,$ \eref{x11} implies 
\beq \label{s11}
s_{11} \, = \, u_{11} \, = \, 0.
\eeq
Next, since the upper left $2$ by $2$ corner of $S$ and $U$ are positive semidefinite, 
\eref{s11} implies that also $s_{12} = u_{12} = 0$ hold, so 
\beq \label{x12}
x_{12} \, = \, t_{12} + v_{12}.
\eeq
Finally, \eref{x12} with  $x_{12} = 0, \, t_{12} \geq 0, \, v_{12} \geq 0$ implies 
$t_{12} = v_{12} = 0, \,$ i.e., $T = V = 0.$ Summarizing, 
$X = S + U \,$ with $S, U \in F_1^*, \, s_{11} = s_{12} = u_{11} = u_{12} = 0, \,$ 
hence $Y = S$ and $Z = U$ are in $H, \,$ as required. 

Next we show that $H$ is a properly minimal face: this comes from the easy-to-check fact that 
$$
F^\perp \, = \, \bpx 0 & 0       & \ti \\
               0 & 0   & \ti \\
             \ti & \ti     & \ti
            \epx.
$$
Finally, we prove that $H$ is not exposed. Let 
$$
Y \, = \, \bpx 0 & 0 & 0 \\ 0 & 1 & 0 \\ 0 & 0 & 0 \epx.
$$
Then clearly $Y \in \ri H, \, $ and 
$$
F \cap Y^\perp \, = \, F \cap H^\perp \, = \, \bpx \oplus & 0 & 0 \\ 0 & 0 & 0 \\ 0 & 0 & 0 \epx.
$$
Choosing $Z$ as the matrix with a $1$ in its upper left corner, and zeros everywhere else, 
we have $Z \in \ri (F \cap H^\perp),$ so the smallest exposed face of $F^*$ that contains $H$ is
$F^* \cap Z^\perp = H + \cone M, \,$ which is strictly larger than $H.$ 

\qed
The following result shows that proving a weaker result suffices to prove Conjecture 
\ref{conj1}. 

\bth \label{weaker-thm}
Suppose that 
\beq \label{k*f}
K^* + F^\perp = F^*
\eeq
holds whenever $K$ is a closed, convex, facially exposed cone, and $F$ is a properly maximal
face of $K.$ Then Conjecture \ref{conj1} is true.
\enth
\pf{} We show that if the assumption of the theorem is true, then \eref{k*f} holds for all $K$ closed, 
convex, facially exposed cones, and all $F$ faces of $K.$ 

Let $K$ be a closed, convex, facially exposed cone. 
We first prove that an $F$ arbitrary face of $K$ is 
facially exposed as a cone in its own right.
Indeed, suppose that $F$ is {\em not} facially exposed. Then there exist 
$F_1$ and $F_2$ faces of $F$ with $F_1 \subsetneq F_2, \,$ and $F^* \cap F_1^\perp = F^* \cap F_2^\perp.$ 
Intersecting both sides of this equation with $K^*$ yields $K^* \cap F_1^\perp = K^* \cap F_2^\perp.$ Since $F_1$ and $F_2$ are  also faces of $K, \,$ this means that $K$ is not facially exposed, a 
contradiction. 

Now, let $F$ again be an arbitrary face of $K.$ 
To show that \eref{k*f} holds for this face, define the chain of faces 
\beq \nonumber
F_0 = K, \, F_1, \, \dots, F_{k-1}, \, F_k = F,
\eeq
where $F_i$ is a properly maximal face of $F_{i-1}$ for $i=1, \dots, k.$ Since all the $F_i$ are facially 
exposed, by the assumption we get 
\beq \nonumber
\ba{rcl}
F_k^* & = & F_{k-1}^* + F_k^\perp, \\
F_{k-1}^* & = & F_{k-2}^* + F_{k-1}^\perp, \\
          & \vdots & \\
F_1^*     & = & F_0^* + F_1^\perp,
\ena
\eeq
hence 
\beq \nonumber
\ba{rcl}
F_k^* & = & F_{k-1}^* + F_k^\perp \\
      & = & F_{k-2}^* + F_{k-1}^\perp + F_k^\perp \\
      & \vdots & \\
      & = & F_0^* + F_1^\perp + \dots + F_k^\perp \\
      & = & F_0^* + F_k^\perp,
\ena
\eeq
as required. 
\qed

\brem {\rm It is known, that if $K$ is a closed, convex cone, and $F$ a properly maximal face of $K, \,$ then 
$F$ is an exposed face of $K$ (\cite[Corollary 2.2]{Tam:85} or \cite[Remark 2.4]{Dick:11}).
We do not use this result, and of course it does not imply that $F$ would be a 
facially exposed cone.}
\erem
If $F$ is a face of $K, \,$ then 
$\lspace F^* = F^\perp, \,$ and if $H$ is a properly minimal face of $F^*, \,$ 
then using part \eref{EE00} in Proposition \ref{lspace-prop}
it follows that $H \cap \lin F$ is an extreme ray of $F^* \cap \lin F.$ 
We define a vector $\xfh$ as the unique vector with norm $1 \,$ that satisfies
$\cone \xfh \, = \, H \cap \lin F \, $ 
(for simplicity, we do not indicate the dependence on $K, \,$ but this should not be confusing).
Then Proposition \ref{lspace-prop} implies 
$$
H \, = \, \cone \xfh + F^\perp.
$$
Also, for an $F$ face of $K$ we denote the orthogonal projection operator onto $\lin F$ by $M_F.$

We first rephrase a condition in Theorem \ref{nice-thm1}. 
\bprop \label{xfhprop}
Let $K$ be a closed, convex cone, $F \tlefteq K, \,$ and $H$ a properly minimal face of $F^*.$ 
Then $K^* \cap H \supsetneq K^* \cap F^\perp$ iff $\xfh \in M_F K^*.$
\eprop
\pf{} We have the following chain of equivalences:
\beq \nonumber
\ba{rclr}
K^* \cap H                       & \supsetneq & K^* \cap F^\perp                        & \LRA \\
(K^* \cap H) \setminus F^\perp   & \neq       & \emp                                    & \LRA \\
(K^* \cap (\cone \xfh + F^\perp)) \setminus F^\perp   & \neq       & \emp             & \LRA \\
\exists \lambda \geq 0, f \in F^\perp: \lambda \xfh + f & \in & K^* \setminus F^\perp        & \LRA \\
\exists \lambda > 0, f \in F^\perp: \lambda \xfh + f & \in & K^* \setminus F^\perp           & \LRA \\
\exists f \in F^\perp:  \xfh + f & \in & K^* \setminus F^\perp                        & \LRA \\
\exists f \in F^\perp:  \xfh + f & \in & K^*                                  & \LRA \\
                               \xfh     & \in & M_F K^*.
\ena
\eeq  
Here the second equivalence comes from \eref{conexfh}, the sixth
from $\xfh \in \lin F \setminus \{ 0 \}, \,$ and the others are trivial. 
\qed

Combining Proposition \ref{xfhprop} with Theorem \ref{nice-thm1} we obtain
\bcor \label{cor-xfh}
Let $K$ be a closed convex cone. Then $K$ is nice, iff 
$\xfh \in M_F K^*$ for all $F \tlefteq K$  and all $H$ properly minimal faces of $F^*.$ 
\ecor
\qed

We now outline a possible avenue of proving Conjecture \ref{conj1}. First, we state
\bprop
Let $K$ be a closed, convex cone, $F \tlefteq K, \,$ and $H$ a properly minimal face of $F^*.$ 
Then the following statements hold.
\benum
\item \label{part1} $\xfh = \lim_i \xfhi, \,$ where 
$H_i$ is a properly minimal, and exposed face of $F^*.$ 
\item \label{part2} If $K$ is 
facially exposed, then for the above $H_i$ we have $\xfhi \in M_F K^*$ for all $i.$
\eenum
\eprop
\pf{} Since $\cone \xfh$ is an extreme ray of $F^* \cap \lin F,$ 
using  Straszewicz's theorem (Theorem 18.6 in \cite{Rockafellar:70}) as it applies to cones 
(see e.g. Theorem 2.12 in \cite{Dick:11}) we get 
\beq \nonumber
\xfh = \lim_i x_i,
\eeq
for some $x_i \in F^* \cap \lin F$ with $\norm{x_i}=1, \,$ and $\cone x_i$ an extreme, exposed ray of 
$F^* \cap \lin F$ for all $i.$ 
By part \eref{EE00} of Proposition \ref{lspace-prop} we get 
$x_i = \xfhi \, $ for some $H_i$ properly minimal faces of $F^*$ for all $i.$ Also, since 
$\cone x_i$ is exposed, by part \eref{again} in Proposition \ref{lspace-prop} 
so is $H_i,$ i.e., $H_i = F^* \cap G_i^\perp $ for some $G_i$ faces of $F$ s.t. $G_i \neq F$ for all $i. \,$ 
This proves part \eref{part1}.

Let us assume that $K$ is facially exposed. Then for the above $G_i$ faces we have 
$K^* \cap G_i^\perp \, \supsetneq \, K^* \cap F^\perp$ for all $i.$ Since 
$K^* \cap H_i \, = \, K^* \cap F^* \cap G_i^\perp \, = \, K^* \cap G_i^\perp, \,$ 
we obtain 
$$
K^* \cap H_i \, \supsetneq \, K^* \cap F^\perp
$$
for all $i.$ Using the equivalence of Lemma \ref{xfhprop}, we get $\xfhi \in M_F K^*$ for all $i,$ 
finishing the proof of \eref{part2}. 
\qed

Summarizing, the facial exposedness of $K$ with the closedness of $M_F K^*$ implies $K^* + F^\perp = F^*.$
Unfortunately, as shown in Proposition 2.1 in \cite{Pataki:07}
the closedness of $M_F K^*$ is actually equivalent to
$K^* + F^\perp = F^*.$ Still, it would be sufficient, and perhaps possible to prove that 
$M_F K^*$ is ``locally'' closed, i.e., when a sequence of vectors from extreme rays of $F^* \cap \lin F$ is in this set, so is their limit.

\nin{\bf Acknowledgement} Thanks are due to Peter J. C. Dickinson for helpful comments on 
copositive and completely positive cones, and to the anonymous 
referees for their thorough reading of the paper, and their helpful comments and suggestions.

\bibliography{D:/bibfiles/mysdp} 


\end{document}